\pdfoutput=1
\documentclass[a4paper,11pt]{amsart}

\usepackage{amsfonts}
\usepackage{amsthm}
\usepackage[colorlinks]{hyperref}
\usepackage{amssymb}
\usepackage{graphicx}
\usepackage{float}
\usepackage[space]{cite}\usepackage{eqparbox,array}
\usepackage{algorithm}
\usepackage[noend]{algorithmic}
\usepackage{dcolumn}
\usepackage{amsmath}
\usepackage{verbatim}
\usepackage{todonotes}%
\usepackage{epstopdf}
\usepackage[utf8]{inputenc}
\usepackage[T1]{fontenc}

\usepackage{tikz}\usetikzlibrary{automata,calc}
\usepackage{epstopdf}
\usepackage{hyphenat}
\usepackage[english]{babel}

\newlength{\lw}
\newlength{\st}
\newlength{\qs}
\newlength{\nd}
\tikzstyle{every picture}+=[auto]
\tikzstyle{every picture}+=[bend angle=10]
\tikzstyle{every picture}+=[join=round]
\tikzstyle{every picture}+=[cap=butt]
\tikzstyle{every picture}+=[line width=\lw]
\tikzstyle{every picture}+=[double distance=2\lw]
\tikzstyle{every picture}+=[shorten >=\st]
\tikzstyle{every picture}+=[node distance=\nd]
\tikzstyle{every loop}=[->,shorten >=\st]
\tikzstyle{place}=[circle,draw,minimum size=\qs]
\tikzstyle{transition}=[rectangle,draw,minimum size=\qs]
\tikzstyle{invisible}=[draw=none,inner sep=0pt,minimum height=0pt]
\newenvironment{petrinet}[1][]{\begin{center}\begin{tikzpicture}}{\end{tikzpicture}\end{center}}

\def\good{/\!\!/}
\newcommand{\PP}{\mathbb{P}}
\newcommand{\CC}{\mathbb{C}}

\begin{document}
\setcounter{footnote}{0}
\setcounter{figure}{0}

\title[Massively parallel computations in algebraic geometry]{Massively parallel computations in algebraic geometry -- not a contradiction}

\author{Janko~B\"ohm}
\address{Janko~B\"ohm, Department of Mathematics, University of Kaiserslautern,  Erwin-Schrödinger-Str., 67663 Kaiserslautern, Germany}
\email{boehm@mathematik.uni-kl.de}

\author{Anne~Fr\"uhbis-Kr\"uger}
\address{Anne~Fr\"uhbis-Kr\"uger, Institut f\"ur algebraische Geometrie, Leibniz Universit\"at Hannover,  Welfengarten 1, 30167 Hannover, Germany}
\email{anne@math.uni-hannover.de}

\author{Mirko~Rahn}
\address{Mirko~Rahn, Competence Center High Performance Computing, Fraunhofer ITWM,  Fraunhofer Platz 1, 67663 Kaiserslautern, Germany}
\email{mirko.rahn@itwm.fhg.de}

 \renewcommand\shortauthors{B\"ohm, J. et al}

\keywords{Computer algebra, Singular, distributed computing, GPI-Space, Petri nets,
computational algebraic geometry, smoothness test, GIT-fan, tropical varieties} 

\thanks{This work has been supported through SPP 1489 and through Project II.5 of SFB-TRR 195 
``Symbolic Tools in Mathematics and their Application'' of the German Research Foundation (DFG)}

\maketitle

\begin{abstract}
The design and implementation of parallel algorithms is a fundamental task in computer algebra. Combining the computer algebra system \textsc{Singular} and the workflow management system GPI-Space, we have developed an infrastructure
for massively parallel computations in commutative algebra and algebraic
geometry. In this note, we give an overview on the current capabilities of this framework by looking into three sample applications: determining smoothness of algebraic varieties, computing GIT-fans in geometric invariant theory, and determining tropicalizations. These applications employ algorithmic methods originating from commutative algebra, sheaf structures on manifolds,  
local geometry, convex geometry, group theory, and combinatorics, illustrating the potential of the framework in further problems in computer algebra.
\end{abstract}

\section{Introduction}

Parallel computations play an increasingly important role in the development of computer infrastructure. In particular, high-performance clusters have the potential for a game-changing
increase in computing power. This raises the question of development and efficient modeling of parallel algorithms in any field relying on large scale calculations. In the context of numerical simulations, massively parallel computations on clusters are nowadays an indispensable tool used for example in weather forecasts or seismic data analysis. In symbolic computing and, in particular, computational algebraic geometry, parallelism also starts to take a center stage \cite{normal, adjoint, farey}. However, due to the unpredictability of time and memory consumption of key algorithmic tools like Buchberger's Algorithm for computing Gr\"obner bases, the use of parallel algorithms has been limited to rather specific contexts. With the goal of a widespread use of parallelism in computer algebra, an effort to adopt the approach of separating the actual computation from the coordination layer has been started. From an approach of this kind, fields like numerical simulation have already benefited significantly in recent years.  Using the computer algebra system {\sc Singular} as the computational back-end within the framework of the workflow management system GPI-space, which employs Petri nets to model the respective algorithm in the coordination layer, the necessary infrastructure has been developed in the PhD thesis of Lukas Ristau, see also \cite{BDFPRR}. So far, we have addressed three sample applications which illustrate the benefits of our approach in computer algebra and will be discussed in this note.

\section{Parallel Computing}\label{sect1}

From the technical point of view, there is more than one paradigm to achieve parallelism. In general, computer scientists identify two different models:  in a shared memory based approach several threads have access to the same data in memory, whereas distributed models run independent processes, possibly on different computers, and communicate their results to the other processes as needed. Not being subject to the obvious restriction of running on the same machine, the latter approach is the suitable one for massively parallel computations, that is for applications scaling up to several hundreds or thousands of parallel computations. On the other hand, it also holds many technical challenges to be met by the expertise of computer scientists, in particular, the need for a clever coordination of the computations to ensure scaling of the performance with the number of cores. 

\section{Framework based on GPI-space and {\sc Singular}}\label{sect2}

For the coordination layer we use the workflow management system GPI-space, which is developed by Fraunhofer ITWM (Kaiserslautern) \cite{GPI}, and allows for automated parallel execution of algorithms. GPI-space provides a scalable runtime system suitable for huge dynamic environments like scientific computing cluster, but works equally well on a sets of computers with heterogeneous hardware, or an individual server. It manages available resources and is tolerant to failure of a node during computation. The virtual memory manager allows the sharing of data; the asynchronous data transfers are managed with the goal of hiding latencies. A key feature is the Petri net based workflow engine, which allows for automatic parallelization and dependency tracking, and will be addressed in the subsequent section.  

As back-end we use {\sc Singular} \cite{Singular}, which is a computer algebra system developed for polynomial computation in commutative algebra, algebraic geometry and singularity theory. Its main workhorse is the Gr\"obner basis engine, around which most of its core algorithms are built. {\sc Singular} exists as a stand-alone software with user interface and as a library version {\tt libSingular}, which is designed for use within other systems like \textsc{sage} or \textsc{OSCAR}. For use in conjunction with GPI-space, we rely on {\tt libSingular} in its existing form, which for this purpose did not have to undergo any significant changes.

\section{Petri nets and parallelism}\label{sect3}

\begin{figure}[h]
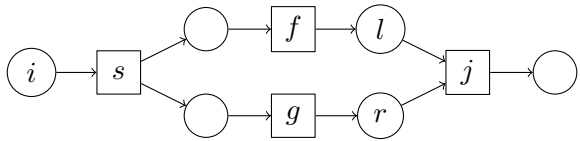

  \begin{petrinet}
  \node[place]      (0)              {$i$};
  \node[transition] (1) [right of=0] {$s$};
  \node[invisible]  (i) [right of=1] {};
  \node[place]      (2) [above of=i,node distance=0.5\nd] {};
  \node[place]      (3) [below of=i,node distance=0.5\nd] {};
  \node[transition] (f) [right of=2] {$f$};
  \node[transition] (g) [right of=3] {$g$};
  \node[place]      (4) [right of=f] {$l$};
  \node[place]      (5) [right of=g] {$r$};
  \node[transition] (j) [right of=i,node distance=3\nd] {$j$};
  \node[place]      (6) [right of=j] {};
  \path[->]
        (0) edge (1)
        (1) edge (2)
        (1) edge (3)
        (2) edge (f)
        (f) edge (4)
        (3) edge (g)
        (g) edge (5)
        (4) edge (j)
        (5) edge (j)
        (j) edge (6)  
  ;
\end{petrinet}
\medskip
  \caption{Example of a Petri net.\label{abb_1}}
\end{figure}

While sequential algorithms are usually described step by step, Petri nets  intrinsically reflect the structure of the algorithms and the state of the computation. As a result they allow for exploiting the parallel structure automatically. A Petri net looks like a directed graph with two kinds of vertices, places (represented as circles in the figures) and transitions (shown as rectangular boxes). The latter contain the elementary functional units, the places can hold marked tokens. These should  be understood to represent pieces of data (in the sense of a so-called \emph{colored} Petri net). Transitions link places, consuming one token from each input place and putting one token onto each output place. So a transition can only fire, if all input places hold a token. Figure~\ref{abb_1} shows a small example of a Petri net. 

A token placed in place $i$ causes transition $s$ to fire which produces one token in each output place. These tokens are then treated independently by transitions 
 $f$ and $g$, which place their output tokens on $l$ and $r$ respectively. The transitions $f$ and $g$ can be executed in parallel, that is, the Petri net allows for task parallelism. The last transition $j$ needs to consume one token from each of $l$ and $r$. Note that the transition $j$ can only fire, if there is at least one token in each of the places $l$ and $r$. If several tokens are available on the places $l$ and $r$ the Petri net does not determine which of them will be consumed first by the transition $j$.  By meeting this rule, we can allow that several parallel instances of each transition run in parallel without destroying the integrity of the computation (data parallelism). By the use of so-called \emph{conditions}, it is possible to ensure that, depending on properties of the tokens, only certain transition can fire.

\section{Applications}
\subsection{Smoothness Test}\label{appl1}

\begin{figure}[h]
  \centering
  \includegraphics[width=0.8\columnwidth]{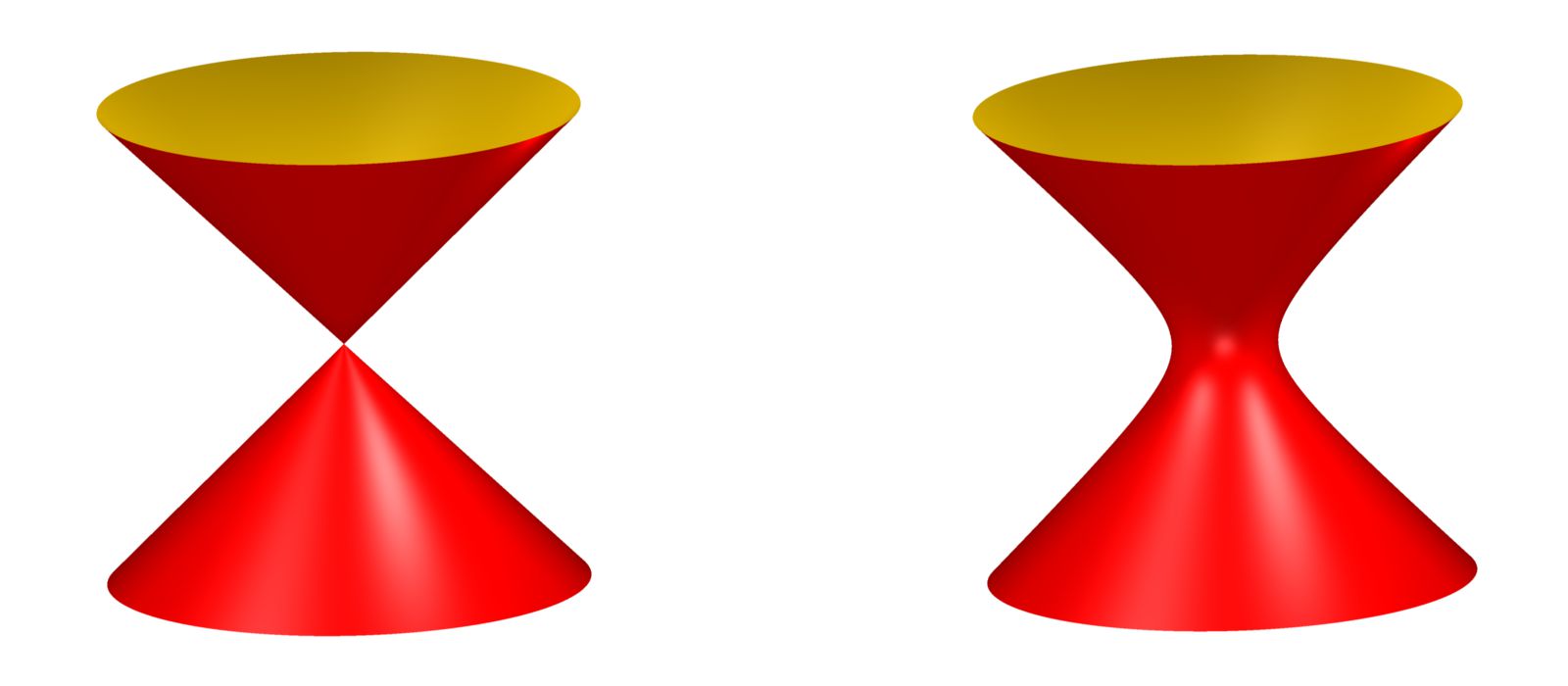}
\medskip
  \caption{A singular and a smooth quadric.\label{abbsing}}
\end{figure}

One of the central tasks of computational algebraic geometry is to explicitly construct objects which exhibit prescribed properties, e.g., general members of moduli spaces or counterexamples to conjectures. In this context, often the question of deciding smoothness of geometric objects arises, since singularities change intrinsic properties of these objects. Thinking in pictures, an algebraic set is smooth, if it looks  in a sufficiently small neighborhood of each of its points like an affine space (that is, like a ${\mathbb C}^k$), see Figure~\ref{abbsing}.

It is classically know that smoothness can be decided by considering the Jacobian ideal generated by all minors of the Jacobian matrix of the correct (that is, codimension) size. Unfortunately, for many problems arising in current research, the Jacobian ideal or even a single minor can be too large to be handled with current computational means; for instance, a recently constructed Godeaux surface (a surface of \emph{general type} with geometric genus and irregularity zero and $K^2=1$) is minimally described by 38 generators in 13 variables, which leads to billions minors.

Taking a completely different approach to deciding smoothness (or in this case more precisely regularity), one can follow ideas of Hironaka's resolution of singularities and use his termination criterion for the resolution process. In a nutshell,  given affine algebraic set $X$, this criterion locally uses a descending induction on the dimension of ambient spaces $W_i\supset X$ and checks at each level and point $p\in W_i$ the value of a certain invariant, the order of the defining ideal of the given affine algebraic set $X$ in ${\mathcal O}_{W_i,p}$. If at some level $i$ and point $p$ the order exceeds one, the point $p$ is singular. From an algorithmic point of view, we can then terminate the computation with a certificate of non-smoothness; if the order at every level and point has the value one until we reach the lowest level where the ideal of $X$ equals the whole ring ${\mathcal O}_{W_k,p}$, we obtain a certificate of smoothness. Transferring this approach from the local analytic setting at individual points to the algebraic setting making use of Zariski open charts, the descending induction requires the computation of a suitable open covering at each level. This leads to a tree of charts with a potentially large number of leaves growing with the level $i$. Moreover, the deeper the level, the more difficult becomes the computation of the open covering and of the current ambient space in each of the new open sets.  

In practice, it is most useful to first proceed with the descending induction in order to pass from a computationally very hard global problem to many computationally simpler local problems. However, before the combinatorial complexity becomes too large, one rather passes in each open set to a relative version of the Jacobian criterion when reaching a given codimension $c$. Since then both the number of generators of the ideal of $X$ and the codimension are smaller, the computation of the ideal of minors (in each of the open set) is a significantly easier task than in the beginning.

\begin{figure}[h]
  \centering
  \includegraphics[width=0.4\columnwidth]{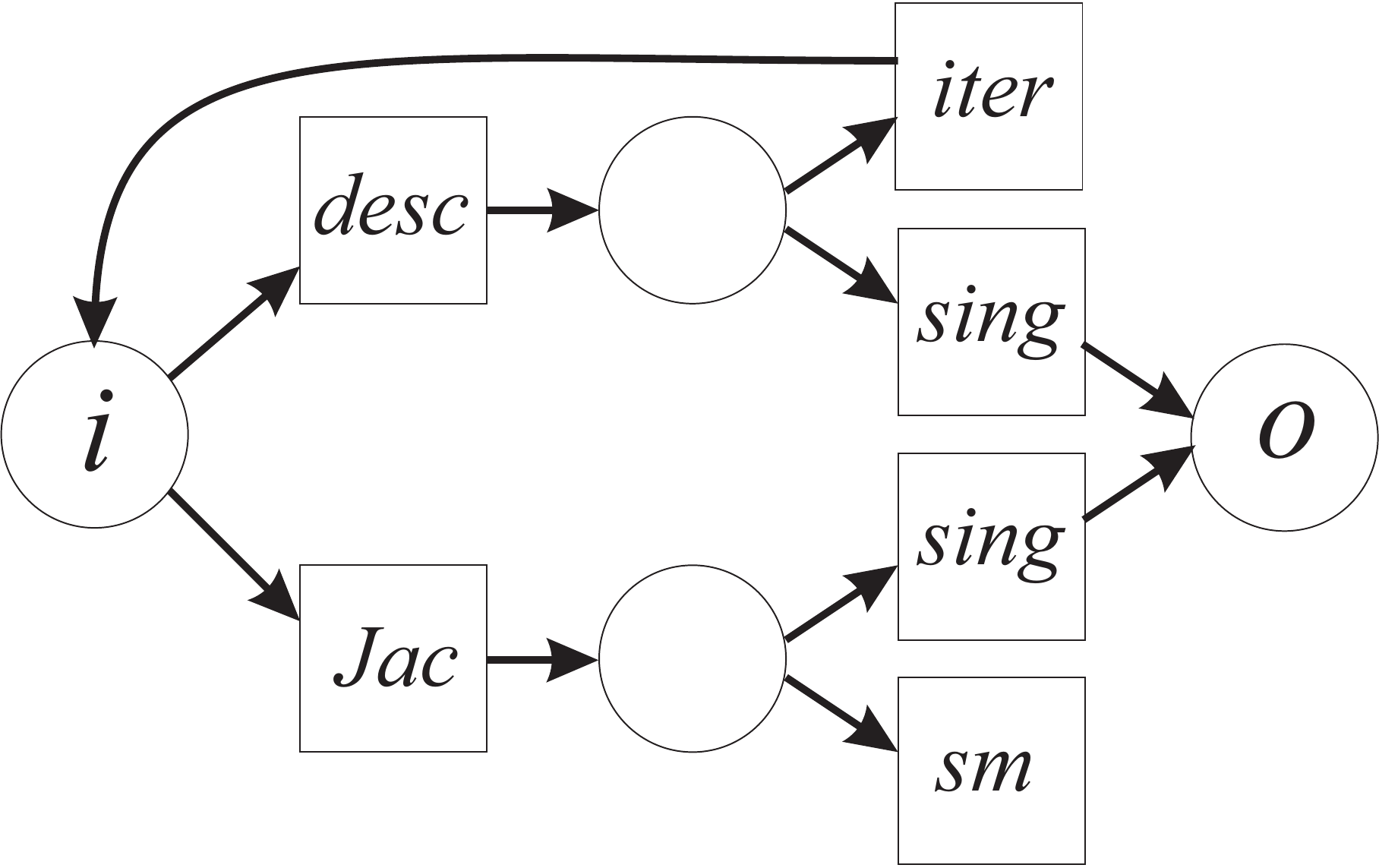}

\medskip
  \caption{Basic structure of the Petri net for the smoothness test.\label{smoothness}}
\end{figure}

A simplified version of the Petri net used for modeling the algorithm is 
shown in Figure \ref{smoothness}. 
Here tokens represent tuples consisting out of $X$, a principal open subset $U$ of affine space and a local ambient algebraic set $W\supset X \cap U$. The input token is placed on \emph{i}. Depending on whether the codimension $c$ has been reached for a given token (which is tested for by using conditions), either  the relative Jacobian criterion (transition \emph{Jac}) or the descent in dimension (transition \emph{desc}) fires. These transitions add to each of their output tokens the knowledge whether a singularity has been detected. If so, the respective transition \emph{sing} fires and places a token on the output place \emph{o}, indicating that the variety is singular. Otherwise the token is deleted by \emph{sm} or replaced on \emph{i}, in the respective cases. If there are no tokens left in the Petri net, we know that the input algebraic set in the input principal open subset was smooth. For details on the algorithm and implementation see \cite{smoothtst, BDFPRR}.

For the Godeaux surface, the global computation terminates due to memory exhaustion on our hardware at $384$ GB of RAM, whereas none of the computations in charts require more than 3.1 GB of RAM and the whole computation takes less than an hour on $128$ cores. 

It is important to note that, due to the massively parallel execution of the Petri net, the computation will automatically determine, from all possible choices of charts, a covering of $X$ which leads to the (non-)smoothness certificate in the fastest possible way. This in fact leads to a super-linear speedup with the number of cores: While the computation on $16$ cores takes about $53000$ seconds, the computation on 
$128$ cores finishes after $3100$ seconds, i.e. using the $8$-fold number of cores we achieve a speedup by a factor of $17$ (with timings measured on a cluster of ITWM with 192 nodes, each with $16$ Intel Xeon E5-2670 cores and $64$ GB of RAM). Note that the maximal speedup is limited by the total number of leaves of the tree of charts produced by the algorithm.

\subsection{Computing GIT-fans}\label{appl2}

\begin{figure}[h]
  \centering
  \includegraphics[width=0.5\columnwidth]{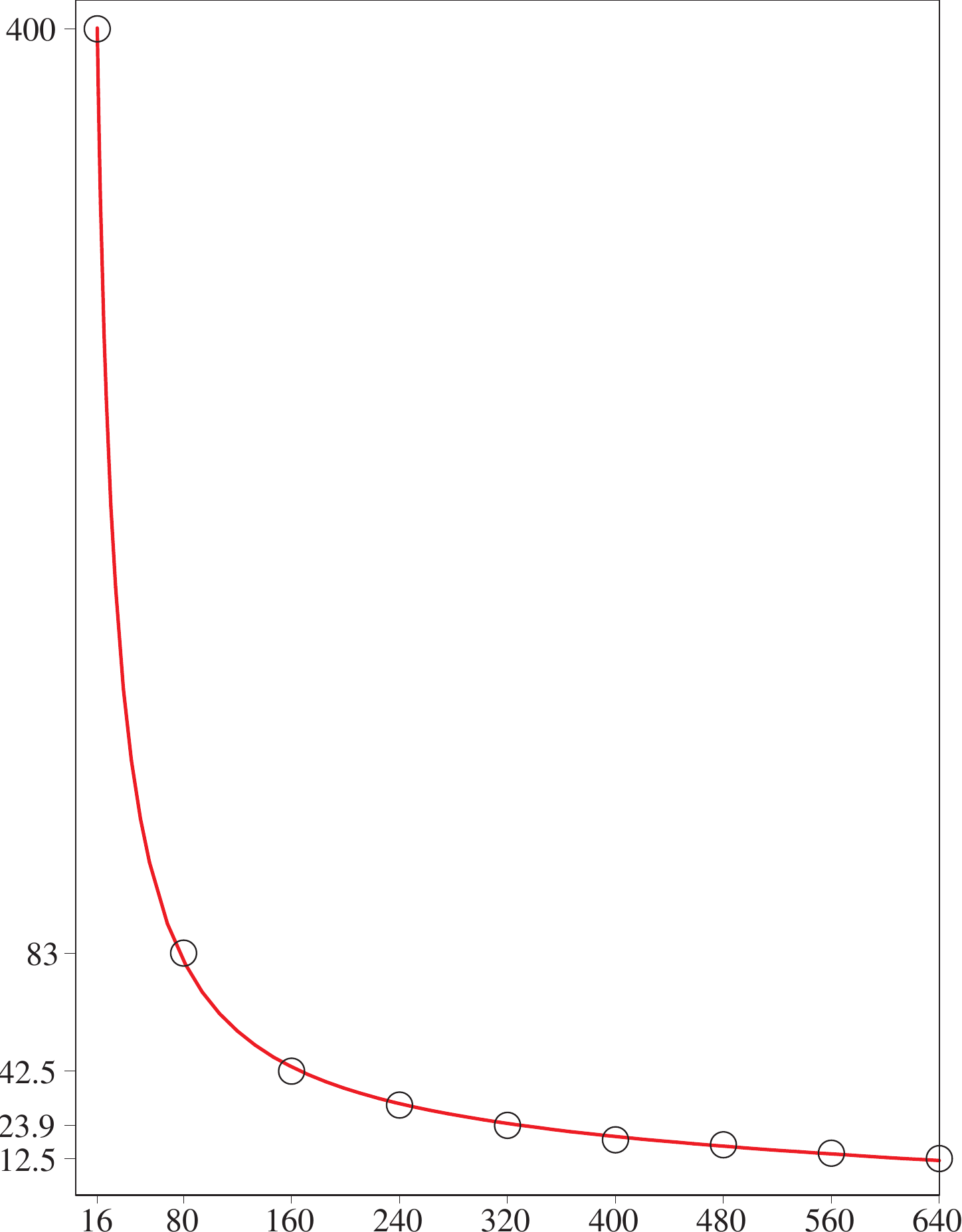}

\bigskip
  \caption{Computation times (in minutes) for the Mori chamber decomposition of $\overline{M}_{0,6}$.\label{GITtimes}}
\end{figure}

\begin{figure}[h]
  \centering
  \includegraphics[width=0.5\columnwidth]{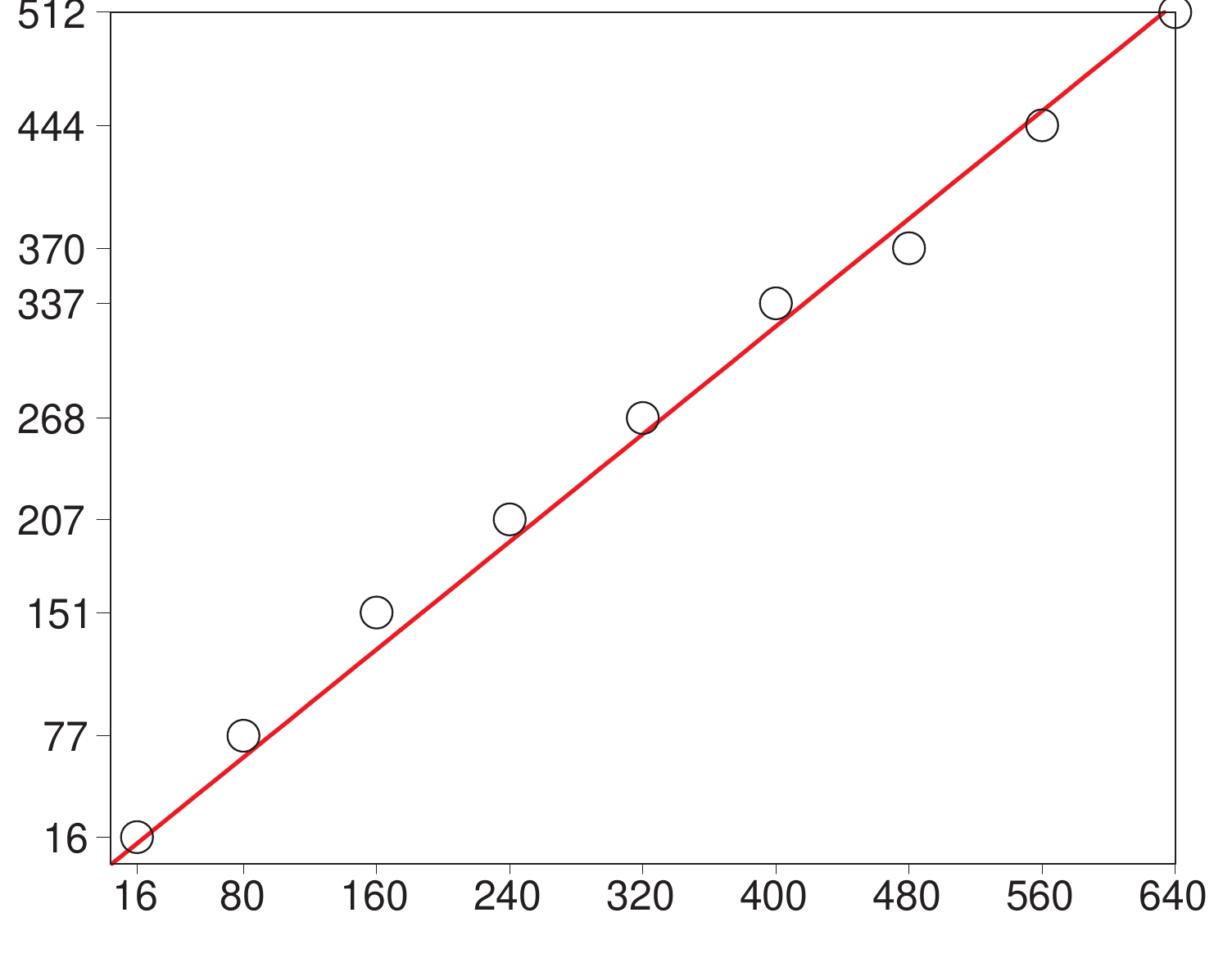}

\bigskip
  \caption{Scaling with the number of cores for computing the Mori chamber decomposition of $\overline{M}_{0,6}$.\label{GITscaling}}
\end{figure}

Geometric Invariant Theory (GIT)  aims at associating a reasonable quotient $X\good G$ to an algebraic variety $X$ on which an algebraic group $G$ acts. This setup occurs, e.g., when a parameterizing space for classes of geometric objects is to be constructed, where the action of the group $G$ on $X$ arises from isomorphisms between the objects. The homogeneous space $X/G$ is not a good candidate for $X\good G$, as it may not be an algebraic variety. For an affine variety $X$, one can rather define the quotient $X\good G$ via the ring of invariant functions of $X$. However, this quotient may show very little structure:
For example, the quotient $\mathbb{C}^2 \good \CC^*$ with the action given by component-wise multiplication is just a point. If we allow for choosing an open subset of $X$, we obtain a much richer geometry, e.g., $(\mathbb{C}^2\setminus \{(0,0)\})\good \CC^* = \PP^1$.
For given $X$, there are usually many choices of open subsets $U\subseteq X$ with different choices leading to
different birationally equivalent quotients
$U\good G$. To describe this behavior, Dolgachev and Hu
\cite{DolgachevHu} introduced the \emph{GIT-fan}, which is a polyhedral fan\footnote{A polyhedral fan is a
finite set of strongly convex rational polyhedral cones such that their faces are again elements of the set and the intersection of any two cones is a face of both.} describing this variation of GIT-quotients. In \cite{BKR} a parallel algorithm for computing GIT-fans, which is also designed to make use of symmetries of the setup, has been described. The algorithm is based on symbolic methods from commutative algebra (Gr\"obner bases), convex geometry (double description) and group theory (orbit decomposition according to the finite symmetry group). Its most important substructure is a traversal of a complete fan which starts at some maximal dimensional cone of the GIT-fan, passes through codimension one faces to all its neighbors, as far as they are not known yet, and iterates. For a short account on the algorithm see \cite{BDKR}. In the Master's thesis of Christian Reinbold \cite{Reinbold}, this algorithm has been modeled in terms of a Petri net and has been implemented using our framework. We have applied this implementation to compute the Mori chamber decomposition of the cone of movable divisors of the Deligne-Mumford moduli space $\overline{M}_{0,6}$ of $6$-pointed stable curves of genus $0$. Figure \ref{GITtimes} shows the computation time plotted against the number of cores in use (on the cluster described above). We observe an impressive linear scaling up to $640$ cores, the maximum number we have tried so far, see Figure~\ref{GITscaling}.


\subsection{Computing tropical varieties}\label{appl3} 

\begin{figure}[h]
  \centering
  \includegraphics[width=0.85\columnwidth]{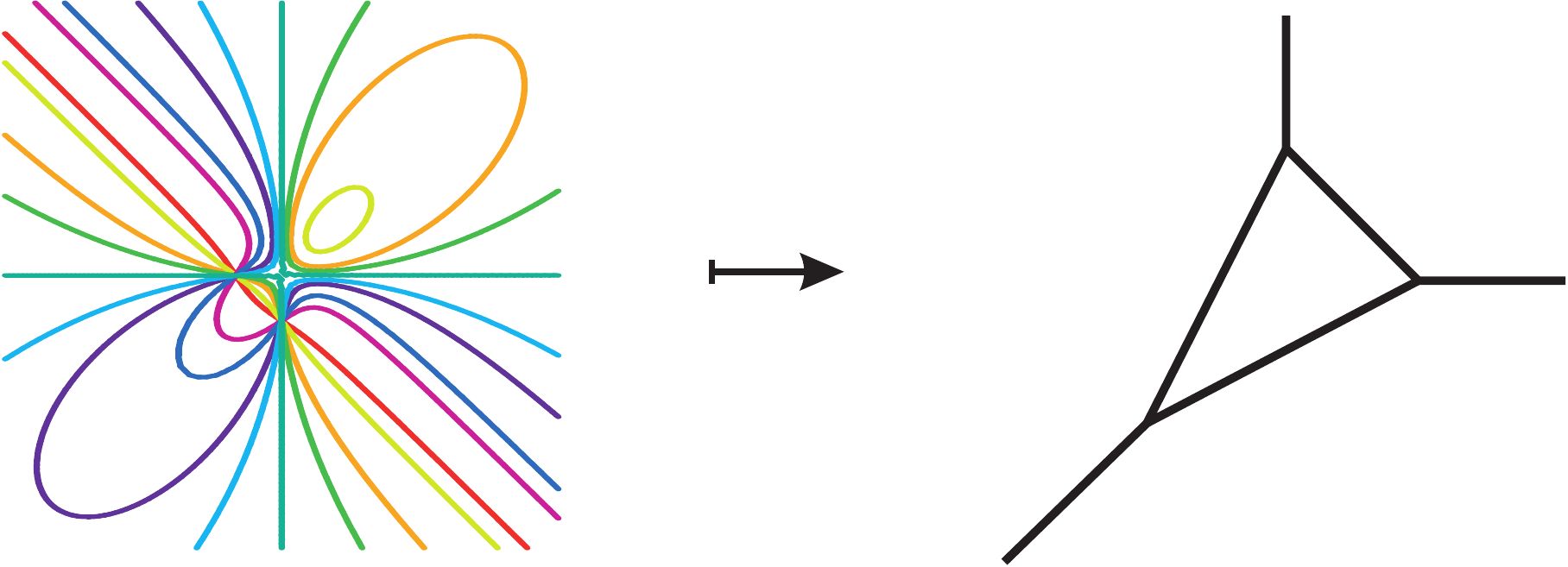}
  \caption{Tropicalization of a family of elliptic curves.\label{tropical}}
\end{figure}

Tropical geometry is a piecewise linear combinatorial version of algebraic geometry. From an element $f$ of a polynomial ring $R=K[x_1,\ldots,x_n]$ over a field $K$ with valuation (usually the Puiseux series field in a variable $t$) one obtains a piecewise linear function by replacing the field operations with those of the tropical semi-ring $$a\oplus b :=\operatorname{min}(a,b),\hspace{7mm} f\odot g := a+b,$$ and coefficients by their valuation. Given an ideal $I\subset R$ and hence an algebraic variety $V(I)$, one can associate to $I$ the tropical variety $T(I)$ defined as the common corner locus of the tropical polynomials corresponding to elements of $I$. By the Bieri-Groves theorem \cite{BG}, the tropical variety corresponds to the set of valuation-tuples of $K$-points of $I$. See Figure \ref{tropical} for a visualization of (fibers in) the family of plane elliptic curves defined by $I=\left\langle t\cdot (x^3+y^3+1) +xy\right\rangle \in K[x,y]$ and the associated tropical variety. 
Computation of tropical varieties is usually a difficult task, both due to the combinatorial complexity of the resulting tropical variety and the computations of the individual faces, which are based on Buchberger's algorithm for finding a Gr\"obner basis. The standard algorithm for this problem starts with one face of the tropical variety and passes to all neighbors by making use of the fact that tropical varieties of prime ideals are irreducible and connected in codimension one \cite{BJSST}. 

Our massively parallel algorithm for computing tropical varieties is based on the fan traversal through codimension one faces (i.e., on a graph expansion problem), see \cite{Bendle}.  Neighboring faces are computed by determining the tropical link via reduction to the curves case and considering Puiseux expansions, see \cite{RH}. The traversal is derived from the fan traversal discussed in Section \ref{appl2}. As in the GIT-fan setting, the algorithm can make use of symmetries.
 However, the algorithm is more involved since, although tropical varieties can also be described as fans, they are not of maximal dimension. As a result, when passing through a codimension one facet, there will typically be more than one neighboring cone of maximal dimension. 

Our implementation computes the previously unknown tropical Grassmannian $\mathcal{G}_{3,8}$ on $768$ cores in less than $20$ minutes with a parallel efficiency comparable to the GIT-fan algorithm. This result opens up the possibility of investigating the relation of the tropical Grassmannian to the Dressian using methods from matroid theory, and the study of the positive Grassmannian.

\section{Conclusion and Outlook}\label{ende}
Using GPI-Space in conjunction with \textsc{Singular} and possibly other computer algebra systems with C-library interface provides an efficient infrastructure for massively parallel computation in computer algebra. By modeling algorithms in terms of Petri nets, parallel structures of the problem can be exploited in a transparent and straight-forward manner. Our infrastructure can make use of the potential of large clusters, but is also suitable for use on smaller multi-core servers and personal computers. Based on the example applications considered so far, we believe that
a multitude of algorithmic problems in computer algebra can benefit
significantly from our approach. Current research addresses, for example, the computation of Hironaka resolutions of singularities \cite{FRS}, Zeta-functions, positive tropical varieties, integration-by-parts identities for Feynman integrals in high-energy physics via the algorithm developed in \cite{F1,F2}, and  generating functions for tropical numerical invariants via the algorithms introduced in \cite{T1,T2}.
%%%%%%%%%%%%%%%%%%%%%%%%%%%%%%%%%%%%%%%%%%%%%%%%%%%%%%%%%%%%%%%%%%%%%%%%%%%%%%%%
% Literature
%%%%%%%%%%%%%%%%%%%%%%%%%%%%%%%%%%%%%%%%%%%%%%%%%%%%%%%%%%%%%%%%%%%%%%%%%%%%%%%%

\end{document}